\documentclass{article}
\usepackage{amssymb}
\usepackage{amsmath}

\newtheorem{theorem}{Theorem}

\newtheorem{definition}[theorem]{Definition}

\newtheorem{lemma}{Lemma}

\newtheorem{remark}{Remark}

\begin{document}

\title{\textbf{On some Leindler's theorem on application of the class }$%
\mathbf{NMCS}$}
\author{\textbf{W\l odzimierz \L enski \ and Bogdan Szal} \\
University of Zielona G\'{o}ra\\
Faculty of Mathematics, Computer Science and Econometrics\\
65-516 Zielona G\'{o}ra, ul. Szafrana 4a\\
P O L A N D\\
W.Lenski@wmie.uz.zgora.pl ,B.Szal @wmie.uz.zgora.pl}
\date{}
\maketitle

\begin{abstract}
We show the results in the class $GM\left( _{5}\beta \right) $ corresponding
to the theorem of L. Leindler [A note on strong approximation of Fourier
series, Analysis Mathematica, 29(2003), 195--199] on strong approximation by
matrix means of Fourier series constructed by the sequences from the class $%
NMCS$.

\ \ \ \ \ \ \ \ \ \ \ \ \ \ \ \ \ \ \ \ 


\ \ \ \ \ \ \ \ \ \ \ \ \ \ \ \ \ \ \ 

\end{abstract}

\section{Introduction}

Let $C_{2\pi }\;$be the class of all $2\pi $--periodic real--valued
functions continuous over $Q=$ $[-\pi ,\pi ]$ with the norm%
\begin{equation*}
\Vert f\Vert :=\sup\limits_{t\in Q}\mid f(t)\mid
\end{equation*}%
and consider the trigonometric Fourier series of $f\in C_{2\pi }$ with the
partial sums\ $S_{k}f$.

Let $A:=\left( a_{n,k}\right) $ be an infinite matrix of real nonnegative
numbers such that%
\begin{equation}
\sum_{k=0}^{\infty }a_{n,k}=1\text{, where }n=0,1,2,...\text{\ ,}  \label{T1}
\end{equation}%
and let the\ $A-$transformation$\ $of \ $\left( S_{k}f\right) $ be given by%
\begin{equation*}
T_{n,A}^{\text{ }}f\left( x\right) :=\sum_{k=0}^{\infty }a_{n,k}S_{k}f\left(
x\right) \text{ \ \ \ }\left( \text{ \ }n=0,1,2,...\right) .
\end{equation*}

Let us consider the strong mean 
\begin{equation*}
T_{n,A}^{p}f\left( x\right) =\sum_{k=0}^{\infty }a_{n,k}\left\vert
S_{k}f\left( x\right) -f\left( x\right) \right\vert ^{p}
\end{equation*}%
and as measures of approximation by such quantity we use the best
approximation of $f$ by trigonometric polynomials $t_{k}$ of order at most $k
$ and the modulus of continuity of $\ f$ defined by the formulas%
\begin{equation*}
E_{k}(f)=\inf_{t_{k}}\left\Vert f-t_{k}\right\Vert 
\end{equation*}%
and%
\begin{equation*}
\omega f\left( \delta \right) =\sup_{\left\vert t\right\vert \leq \delta
}\left\Vert f\left( \cdot +t\right) -f\left( \cdot \right) \right\Vert ,
\end{equation*}%
respectively.

In \cite{MT} S. M. Mazhar and V. Totik proved the following theorem:

\begin{theorem}
Suppose $A:=\left( a_{n,k}\right) $ satisfies (\ref{T1}), $lim_{n\rightarrow
\infty }a_{n,0}=0$ and%
\begin{equation*}
a_{n,k}\geq a_{n,k+1}\text{ \ \ }k=0,1,2,...\text{ \ \ }n=0,1,2,...,
\end{equation*}%
then%
\begin{equation*}
\left\Vert T_{n,A}^{\text{ }}f-f\right\Vert \leq K\sum_{k=0}^{\infty
}a_{n,k}\omega f\left( \frac{1}{k+1}\right) .
\end{equation*}
\end{theorem}

Recently, L. Leindler \cite{3} defined a new class of sequences named as
sequences of rest bounded variation, briefly denoted by $RBVS$, i.e.,%
\begin{equation}
RBVS=\left\{ a:=\left( a_{n}\right) \in 
\mathbb{C}
:\sum\limits_{k=m}^{\infty }\left\vert a_{k}-a_{k+1}\right\vert \leq
K\left( a\right) \left\vert a_{m}\right\vert \text{ for all }m\in 
\mathbb{N}
\right\} ,  \label{1}
\end{equation}%
where here and throughout the paper $K\left( a\right) $ always indicates a
constant only depending on $a$.

Denote by $MS$ the class of monotone decreasing sequences and $CQMS$ the
class of classic quasimonotone decreasing sequences ($a\in CQMS$ means that $%
\left( a_{n}\right) \in 
\mathbb{R}
_{+}$ and there exists an $\alpha >0$ such that $a_{n}/n^{\alpha }$ is
decreasing), then it is obvious that%
\begin{equation*}
MS\subset RBVS\cap CQMS.
\end{equation*}%
L. Leindler \cite{16} proved that the class $CQMS$ and $RBVS$ are not
comparable. In \cite{8} L. Leindler considered the class of mean rest
bounded variation sequences $MRBVS$, where%
\begin{equation*}
MRBVS=\left\{ a:=\left( a_{n}\right) \in 
\mathbb{C}
:\right.
\end{equation*}%
\begin{equation}
\left. \sum\limits_{k=m}^{\infty }\left\vert a_{k}-a_{k+1}\right\vert \leq
K\left( a\right) \frac{1}{m}\sum\limits_{k\geq m/2}^{m}\left\vert
a_{k}\right\vert \text{ for all }m\in 
\mathbb{N}
\right\} .  \label{2}
\end{equation}%
It is clear that%
\begin{equation*}
RBVS\subseteq MRBVS.
\end{equation*}%
In \cite{S} the second author proved that $RBVS\neq MRBVS$. Moreover, the
above theorem was generalized for the class $MRBVS$ in \cite{Sz} .

Further, the class of general monotone coefficients, $GM$, is defined as
follows ( see \cite{10}):%
\begin{equation}
GM=\left\{ a:=\left( a_{n}\right) \in 
\mathbb{C}
:\sum\limits_{k=m}^{2m-1}\left\vert a_{k}-a_{k+1}\right\vert \leq K\left(
a\right) \left\vert a_{m}\right\vert \text{ for all }m\in 
\mathbb{N}
\right\} .  \label{3}
\end{equation}%
It is clear%
\begin{equation*}
RBVS\cup CQMS\subset GM\text{.}
\end{equation*}%
In \cite{6, 10, 11, 12} was defined the class of $\beta -$general monotone
sequences as follows:

\begin{definition}
Let $\beta :=\left( \beta _{n}\right) $ be a nonnegative sequence. The
sequence of complex numbers $a:=\left( a_{n}\right) $ is said to be $\beta -$%
general monotone, or $a\in GM\left( \beta \right) $, if the relation%
\begin{equation}
\sum\limits_{k=m}^{2m-1}\left\vert a_{k}-a_{k+1}\right\vert \leq K\left(
a\right) \beta _{m}  \label{4}
\end{equation}%
holds for all $m$.
\end{definition}

In the paper \cite{12} Tikhonov considered the following examples of the
sequences $\beta _{n}:$

(1) $_{1}\beta _{n}=\left\vert a_{n}\right\vert ,$

(2) $_{2}\beta _{n}=\sum\limits_{k=n}^{n+N}\left\vert a_{k}\right\vert $
for some integer $N$,

(3) $_{3}\beta _{n}=\sum\limits_{\nu =0}^{N}\left\vert a_{c^{\nu
}n}\right\vert $ for some integers $N$ and $c>1$,

(4) $_{4}\beta _{n}=\left\vert a_{n}\right\vert +\sum\limits_{k=n+1}^{\left[
cn\right] }\frac{\left\vert a_{k}\right\vert }{k}$ for some $c>1$,

(5) $_{5}\beta _{n}=\sum\limits_{k=\left[ n/c\right] }^{\left[ cn\right] }%
\frac{\left\vert a_{k}\right\vert }{k}$ for some $c>1$.

It is clear that $GM\left( _{1}\beta \right) =GM.$ Moreover (see \cite[%
Remark 2.1]{12})%
\begin{equation*}
GM\left( _{1}\beta +_{2}\beta +_{3}\beta +_{4}\beta +_{5}\beta \right)
\equiv GM\left( _{5}\beta \right) .
\end{equation*}

Consequently, we assume that the sequence $\left( K\left( \alpha _{n}\right)
\right) _{n=0}^{\infty }$ is bounded, that is, that there exists a constant $%
K$ such that%
\begin{equation*}
0\leq K\left( \alpha _{n}\right) \leq K
\end{equation*}%
holds for all $n$, where $K\left( \alpha _{n}\right) $ denote the sequence
of constants appearing in the inequalities (\ref{1})-(\ref{4}) for the
sequences $\alpha _{n}:=\left( a_{nk}\right) _{k=0}^{\infty }$.

Now we can give the conditions to be used later on. We assume that for all $%
n $%
\begin{equation}
\sum\limits_{k=m}^{2m-1}\left\vert a_{n,k}-a_{n,k+1}\right\vert \leq
K\sum\limits_{k=[m/c]}^{[cm]}\frac{a_{n,k}}{k}  \label{5}
\end{equation}%
holds if $\alpha _{n}=\left( a_{n,k}\right) _{k=0}^{\infty }$ belongs to $%
GM\left( _{5}\beta \right) $, for $n=1,2,...$

Following by L. Leindler \cite{L} a sequence $a:=\left( a_{n}\right) $ of
nonnegative numbers is called a Nearly Monotone Convergent Sequence, or
briefly $a\in NMCS$, if 
\begin{equation*}
\sum\limits_{k=1}^{\infty }a_{k}<\infty \text{ \ \ and \ \ }%
ka_{k}\rightarrow 0\text{ \ as \ }k\rightarrow \infty ,
\end{equation*}%
for all positive integer $r.$

\begin{remark}
If $\sum\limits_{k=1}^{\infty }\left\vert a_{k}\right\vert <\infty $ and $%
\left( a_{k}\right) \in GM\left( _{5}\beta \right) $ then $\left(
a_{k}\right) \in NMCS$.
\end{remark}

The deviation $H_{n,A}^{p}f$ was estimated by L. Leindler in \cite{L} as
follows:

\begin{theorem}
\cite{L} If $f\in C_{2\pi }$, $p>0$, $\left( a_{n,k}\right) _{k=0}^{\infty
}\in NMCS$ for all $n$, and $lim_{n\rightarrow \infty }a_{n,0}=0$ holds,
then 
\begin{equation*}
\left\Vert T_{n,A}^{p}f\right\Vert \leq K\sum_{k=0}^{\infty
}a_{n,k}E_{k}^{p}(f).
\end{equation*}
\end{theorem}

In this note we show that the class $NMCS$ is not proper for the above
estimate.

In our theorem we consider the class $GM\left( _{5}\beta \right) $ instead
of $NMCS.$ Thus we essentially extend the result of S. M. Mazhar and V.
Totik (see \cite{MT}).

We shall write $I_{1}\ll I_{2}$ if there exists a positive constant $K$ ,
sometimes depended on some parameters, such that $I_{1}\leq KI_{2}$.

\section{Statement of the results}

Our main result is the following

\begin{theorem}
If $f\in C_{2\pi }$, $p>0$, $\left( a_{n,k}\right) _{k=0}^{\infty }\in
GM\left( _{5}\beta \right) $ for all $n$, (\ref{T1}) and $lim_{n\rightarrow
\infty }a_{n,0}=0$ hold, then%
\begin{equation}
\left\Vert T_{n,A}^{p}f\right\Vert \ll \sum_{k=0}^{\infty }a_{n,k}E_{\left[ 
\frac{k}{2^{\left[ c\right] }}\right] }^{p}(f)  \label{6}
\end{equation}%
for some $c>1$.
\end{theorem}

\begin{remark}
If we suppose that $\left( a_{n,k}\right) _{k=0}^{\infty }\in MS$ then from (%
\ref{6}) we deduce%
\begin{equation*}
\left\Vert T_{n,A}^{p}f\right\Vert \ll \sum_{k=0}^{\infty
}a_{n,k}E_{k}^{p}(f).
\end{equation*}
\end{remark}

Using the Jackson Theorem \cite[Theorem 13.6]{Z} we can obtain the following
remark.

\begin{remark}
Under the assumptions of Theorem 3%
\begin{equation}
\left\Vert T_{n,A}^{p}f\right\Vert \ll \sum_{k=0}^{\infty }a_{n,k}\omega
^{p}f\left( \frac{\pi }{k+1}\right) .  \label{7}
\end{equation}
\end{remark}

\begin{remark}
We can observe that taking $a_{nn}=1$ and $a_{n,k}=0$ for $k\neq n$ we have $%
\left( a_{n,k}\right) _{k=0}^{n}\in NMCS$ but thus, by Theorem 2\textbf{,}
we obtain the estimate 
\begin{equation*}
\left\Vert S_{n}f-f\right\Vert \ll E_{n}(f)
\end{equation*}%
which is not true in general.
\end{remark}

\begin{remark}
By the considerations similar to these in \cite{L} we can obtain the
estimates 
\begin{equation*}
\left\Vert T_{n,A}^{\varphi }f\right\Vert \ll \sum_{k=0}^{\infty
}a_{n,k}\varphi \left( E_{\left[ \frac{k}{2^{\left[ c\right] }}\right]
}(f)\right) 
\end{equation*}%
and%
\begin{equation*}
\left\Vert T_{n,A}^{\varphi }f\right\Vert \ll \sum_{k=0}^{\infty
}a_{n,k}\varphi \left( \omega f\left( \frac{\pi }{k+1}\right) \right) 
\end{equation*}%
instead of (\ref{6}) and (\ref{7}) respectively, where%
\begin{equation*}
T_{n,A}^{\varphi }f\left( x\right) =\sum_{k=0}^{\infty }a_{n,k}\varphi
\left( \left\vert S_{k}f\left( x\right) -f\left( x\right) \right\vert
\right) 
\end{equation*}
with a nonnegative monotone increasing continuous function $\varphi (t)$ $%
(t\in \lbrack 0,\infty ))$ satisfying the conditions%
\begin{equation*}
\varphi (0)=0,\varphi (t)\leq e^{At},\text{ \ \ }t\in (0,\infty )
\end{equation*}%
and%
\begin{equation*}
\varphi (2t)\leq A\varphi (t),\text{ \ \ }t\in (0,1),
\end{equation*}%
with some constant $A.$
\end{remark}

\section{Auxiliary result}

We shall use the following

\begin{lemma}
(see [\cite{LL}, Theorem 1.11 ). Suppose that $n=O\left( \lambda _{n}\right)
.$ Then, for any continuous function $f$ and for any number $p>0,$we have%
\begin{equation*}
\left\Vert \left\{ \frac{1}{\lambda _{n}}\sum_{k=n-\lambda
_{n}}^{n-1}\left\vert S_{k}f-f\right\vert ^{p}\right\} ^{1/p}\right\Vert \ll
E_{n-\lambda _{n}}(f).
\end{equation*}
\end{lemma}

\section{Proofs of the results}

\subsection{Proof of Theorem 3}

Let%
\begin{equation*}
\left\Vert T_{n,A}^{p}f\right\Vert =\left\Vert \sum_{k=0}^{2^{\left[ c\right]
}-1}a_{n,k}\left\vert S_{k}f-f\right\vert ^{p}+\sum_{k=2^{[c]}}^{\infty
}a_{n,k}\left\vert S_{k}f-f\right\vert ^{p}\right\Vert 
\end{equation*}%
\begin{equation*}
\leq \left\Vert \sum_{k=0}^{2^{\left[ c\right] }-1}a_{n,k}\left\vert
S_{k}f-f\right\vert ^{p}\right\Vert +\left\Vert \sum_{m=\left[ c\right]
}^{\infty }\sum_{k=2^{m}}^{2^{m+1}-1}a_{n,k}\left\vert S_{k}f-f\right\vert
^{p}\right\Vert =I_{1}+I_{2}.
\end{equation*}%
for some $c>1$. Using Lemma we get%
\begin{eqnarray*}
I_{1} &\leq &\left\Vert \sum_{k=0}^{2^{\left[ c\right] }-1}a_{n,k}\frac{k-%
\left[ k/2^{\left[ c\right] }\right] +1}{k-\left[ k/2^{\left[ c\right] }%
\right] +1}\sum\limits_{l=\left[ k/2^{\left[ c\right] }\right]
}^{k}\left\vert S_{l}f-f\right\vert ^{p}\right\Vert  \\
&\leq &2^{\left[ c\right] }\left\Vert \sum_{k=0}^{2^{\left[ c\right]
}-1}a_{n,k}\frac{1}{k-\left[ k/2^{\left[ c\right] }\right] +1}\sum\limits_{l=%
\left[ k/2^{\left[ c\right] }\right] }^{k}\left\vert S_{l}f-f\right\vert
^{p}\right\Vert  \\
&\ll &\sum_{k=0}^{2^{\left[ c\right] }-1}a_{n,k}E_{\left[ k/2^{\left[ c%
\right] }\right] }^{p}(f).
\end{eqnarray*}%
By partial summation, our Lemma gives%
\begin{eqnarray*}
I_{2} &=&\left\Vert \sum_{m=[c]}^{\infty }\left[ \sum_{k=2^{m}}^{2^{m+1}-2}%
\left( a_{n,k}-a_{n,k+1}\right) \sum_{l=2^{m}}^{k}\left\vert
S_{l}f-f\right\vert ^{p}\right. \right.  \\
&&\left. \left. +a_{n,2^{m+1}-1}\sum_{l=2^{m}}^{2^{m+1}-1}\left\vert
S_{l}f-f\right\vert ^{p}\right] \right\Vert 
\end{eqnarray*}%
\begin{eqnarray*}
&\ll &\sum_{m=[c]}^{\infty }\left[ 2^{m}\sum_{k=2^{m}}^{2^{m+1}-2}\left\vert
a_{n,k}-a_{n,k+1}\right\vert E_{2^{m}}^{p}(f)\right.  \\
&&\left. +2^{m}a_{n,2^{m+1}-1}E_{2^{m}}^{p}(f)\right]  \\
&\ll &\sum_{m=[c]}^{\infty }2^{m}E_{2^{m}}^{p}(f)\left[
\sum_{k=2^{m}}^{2^{m+1}-2}\left\vert a_{n,k}-a_{n,k+1}\right\vert
+a_{n,2^{m+1}-1}\right] .
\end{eqnarray*}%
Since (\ref{5}) holds, we have%
\begin{eqnarray*}
&&a_{n,s+1}-a_{n,r} \\
&\leq &\left\vert a_{n,r}-a_{n,s+1}\right\vert \leq \sum_{k=r}^{s}\left\vert
a_{n,k}-a_{n,k+1}\right\vert  \\
&\leq &\sum_{k=2^{m}}^{2^{m+1}-2}\left\vert a_{n,k}-a_{n,k+1}\right\vert \ll
\sum\limits_{k=[2^{m}/c]}^{[c2^{m}]}\frac{a_{n,k}}{k}\text{ \ \ }\left(
2\leq 2^{m}\leq r\leq s\leq 2^{m+1}-2\right) ,
\end{eqnarray*}%
whence%
\begin{equation*}
a_{n,s+1}\ll a_{n,r}+\sum\limits_{k=[2^{m}/c]}^{[c2^{m}]}\frac{a_{n,k}}{k}%
\text{ \ }\left( 2\leq 2^{m}\leq r\leq s\leq 2^{m+1}-2\right) 
\end{equation*}%
and%
\begin{eqnarray*}
2^{m}a_{n,2^{m+1}-1} &=&\frac{2^{m}}{2^{m}-1}%
\sum_{r=2^{m}}^{2^{m+1}-2}a_{n,2^{m+1}-1} \\
&\ll &\sum_{r=2^{m}}^{2^{m+1}-2}\left(
a_{n,r}+\sum\limits_{k=[2^{m}/c]}^{[c2^{m}]}\frac{a_{n,k}}{k}\right)  \\
&\ll
&\sum_{r=2^{m}}^{2^{m+1}-1}a_{n,r}+2^{m}\sum\limits_{k=[2^{m}/c]}^{[c2^{m}]}%
\frac{a_{n,k}}{k}
\end{eqnarray*}%
whence%
\begin{equation*}
I_{2}\ll \sum_{m=[c]}^{\infty }\left\{
2^{m}E_{2^{m}}^{p}(f)\sum\limits_{k=[2^{m}/c]}^{[c2^{m}]}\frac{a_{n,k}}{k}%
+E_{2^{m}}^{p}(f)\sum_{k=2^{m}}^{2^{m+1}-1}a_{n,k}\right\} .
\end{equation*}

Finally, by elementary calculations we get%
\begin{eqnarray*}
I_{2} &\ll &\sum_{m=[c]}^{\infty }\left\{
2^{m}E_{2^{m}}^{p}(f)\sum\limits_{k=2^{m-\left[ c\right] }}^{2^{m+\left[ c%
\right] }}\frac{a_{n,k}}{k}+E_{2^{m}}^{p}(f)\sum_{k=2^{m}}^{2^{m+1}}a_{n,k}%
\right\} \\
&\ll &\sum_{m=[c]}^{\infty }E_{2^{m}}^{p}(f)\sum\limits_{k=2^{m-\left[ c%
\right] }}^{2^{m+\left[ c\right] }}a_{n,k} \\
&=&\sum_{m=[c]}^{\infty }E_{2^{m}}^{p}(f)\sum\limits_{k=2^{m-\left[ c\right]
}}^{2^{m}-1}a_{n,k}+\sum_{m=[c]}^{\infty
}E_{2^{m}}^{p}(f)\sum\limits_{k=2^{m}}^{2^{m+\left[ c\right] }}a_{n,k}
\end{eqnarray*}%
\begin{equation*}
\ll \sum_{m=[c]}^{\infty }\sum\limits_{k=2^{m-\left[ c\right]
}}^{2^{m}-1}a_{n,k}E_{k}^{p}(f)+\sum_{m=[c]}^{\infty
}\sum\limits_{k=2^{m}}^{2^{m+\left[ c\right] }}a_{n,k}E_{\left[ \frac{k}{2^{%
\left[ c\right] }}\right] }^{p}(f)
\end{equation*}%
\begin{equation*}
\ll \sum_{m=[c]}^{\infty }\sum\limits_{k=2^{m-\left[ c\right]
}}^{2^{m}-1}a_{n,k}E_{k}^{p}(f)+\sum_{m=[c]}^{\infty
}\sum\limits_{k=2^{m}}^{2^{m+\left[ c\right] }-1}a_{n,k}E_{\left[ \frac{k}{%
2^{\left[ c\right] }}\right] }^{p}(f)+\sum_{m=[c]}^{\infty
}E_{2^{m}}^{p}(f)a_{n,2^{m+\left[ c\right] }}
\end{equation*}%
\begin{eqnarray*}
&=&\sum_{m=[c]}^{\infty }\sum\limits_{r=1}^{\left[ c\right]
}\sum\limits_{k=2^{m-r}}^{2^{m-r+1}-1}a_{n,k}E_{k}^{p}(f)+\sum_{m=[c]}^{%
\infty }\sum\limits_{r=0}^{\left[ c\right] -1}\sum%
\limits_{k=2^{m+r}}^{2^{m+r+1}-1}a_{n,k}E_{\left[ \frac{k}{2^{\left[ c\right]
}}\right] }^{p}(f) \\
&&+\sum_{m=[c]}^{\infty }E_{2^{m}}^{p}(f)a_{n,2^{m+\left[ c\right] }}
\end{eqnarray*}%
\begin{eqnarray*}
&\leq &\sum\limits_{r=1}^{\left[ c\right] }\sum\limits_{k=2^{\left[ c%
\right] -r}}^{\infty }a_{n,k}E_{k}^{p}(f)+\sum\limits_{r=0}^{\left[ c\right]
-1}\sum\limits_{k=2^{\left[ c\right] +r}}^{\infty }a_{n,k}E_{\left[ \frac{k%
}{2^{\left[ c\right] }}\right] }^{p}(f)+\sum\limits_{k=2^{2\left[ c\right]
}}^{\infty }a_{n,k}E_{\left[ \frac{k}{2^{\left[ c\right] }}\right] }^{p}(f)
\\
&\ll &\sum\limits_{k=0}^{\infty }a_{n,k}E_{\left[ \frac{k}{2^{\left[ c%
\right] }}\right] }^{p}(f)
\end{eqnarray*}

Thus we obtain the desired result. $\square $

\subsection{Proof of Remark 1}

For $j=k+2,k+3,...,2k$ we get%
\begin{equation*}
\sum\limits_{l=k}^{j-2}\left\vert a_{l}-a_{l+1}\right\vert \geq
\sum\limits_{l=k}^{j-2}\left\vert \left\vert a_{l}\right\vert -\left\vert
a_{l+1}\right\vert \right\vert \geq \sum\limits_{l=k}^{j-2}\left\vert
a_{l}\right\vert -\sum\limits_{l=k+1}^{j-1}\left\vert a_{l}\right\vert
=\left\vert a_{k}\right\vert -\left\vert a_{j-1}\right\vert
\end{equation*}%
Therefore, if $\left( a_{k}\right) \in GM\left( _{5}\beta \right) $ then%
\begin{equation*}
\left\vert a_{k}\right\vert \leq \sum\limits_{l=k}^{j-2}\left\vert
a_{l}-a_{l+1}\right\vert +\left\vert a_{j-1}\right\vert
\end{equation*}%
\begin{equation*}
\leq \sum\limits_{l=\left[ j/2\right] }^{2\left[ j/2\right] -1}\left\vert
a_{l}-a_{l+1}\right\vert +\left\vert a_{j-1}\right\vert
\end{equation*}%
\begin{equation*}
\ll \sum\limits_{l=\left[ j/2c\right] }^{\left[ cj/2\right] }\frac{%
\left\vert a_{l}\right\vert }{l}+\left\vert a_{j-1}\right\vert \ll \frac{1}{k%
}\sum\limits_{l=\left[ k/2c\right] }^{\left[ ck\right] }\left\vert
a_{l}\right\vert +\left\vert a_{j-1}\right\vert .
\end{equation*}%
Summing up on $j$ and using the assumption $\sum\limits_{k=1}^{\infty
}\left\vert a_{k}\right\vert <\infty $ we get for $k>1$ that%
\begin{equation*}
k\left\vert a_{k}\right\vert =\frac{k}{k-1}\sum\limits_{j=k+2}^{2k}\left%
\vert a_{k}\right\vert \ll \sum\limits_{j=k+2}^{2k}\left( \frac{1}{k}%
\sum\limits_{l=\left[ k/2c\right] }^{\left[ ck\right] }\left\vert
a_{l}\right\vert +\left\vert a_{j-1}\right\vert \right)
\end{equation*}%
\begin{equation*}
\ll \sum\limits_{l=\left[ k/2c\right] }^{\left[ ck\right] }\left\vert
a_{l}\right\vert +\sum\limits_{j=k+2}^{2k}\left\vert a_{j-1}\right\vert
=\sum\limits_{l=\left[ k/2c\right] }^{\left[ ck\right] }\left\vert
a_{l}\right\vert +\sum\limits_{j=k+1}^{2k}\left\vert a_{j}\right\vert
\end{equation*}%
\begin{equation*}
\leq 2\sum\limits_{l=\left[ k/2c\right] }^{\left[ ck\right] +2k}\left\vert
a_{l}\right\vert \rightarrow 0,
\end{equation*}%
whence $\left( a_{k}\right) \in NMCS.$ $\square $

\subsection{Proof of Remark 2}

If $\left( a_{n,k}\right) _{k=0}^{\infty }\in MS$ then $\left(
a_{n,k}\right) _{k=0}^{\infty }\in GM\left( _{5}\beta \right) $ and using
Theorem 3 we obtain%
\begin{eqnarray*}
\left\Vert T_{n,A}^{p}f\right\Vert  &\ll &\sum_{k=0}^{\infty }a_{n,k}E_{
\left[ \frac{k}{2^{\left[ c\right] }}\right] }^{p}(f)=\sum_{k=0}^{\infty
}\sum\limits_{m=k2^{\left[ c\right] }}^{\left( k+1\right) 2^{\left[ c\right]
}-1}a_{n,m}E_{\left[ \frac{m}{2^{\left[ c\right] }}\right] }^{p}(f) \\
&=&\sum_{k=0}^{\infty }E_{k}^{p}(f)\sum\limits_{m=k2^{\left[ c\right]
}}^{\left( k+1\right) 2^{\left[ c\right] }-1}a_{n,m}\leq \sum_{k=0}^{\infty
}2^{\left[ c\right] }E_{k}^{p}(f)a_{n,k2^{\left[ c\right] }} \\
&\leq &\left\{ 2^{\left[ c\right] }\right\} ^{1/p}\sum_{k=0}^{\infty
}E_{k}^{p}(f)a_{n,k}.
\end{eqnarray*}%
This ends our proof. $\square $

\subsection{Proof of Remark 5}

The proof is similar to the proof of Theorem 3. The difference is such that
we use the following Totik estimate (see \cite{T})%
\begin{equation*}
\frac{1}{n}\sum_{k=n+1}^{2n}\varphi \left( \left\vert S_{k}f\left( x\right)
-f\left( x\right) \right\vert \right) \leq K\varphi \left( E_{n}\left(
f\right) \right) 
\end{equation*}%
instead of the inequality from Lemma. $\square $

\end{document}